\documentclass[11pt,a4paper,language=english]{bilingualmath2}

\titleen{Generalized Cluster Algebra Homomorphisms Induced by Cluster Tube Embeddings}
\paperauthor{Mingze Li}
\affiliationen{Beijing Normal University}
\paperemail{mathlmz@mail.bnu.edu.cn}
\dateen{\englishtoday}

\newcommand{\add}{\operatorname{add}}
\newcommand{\ind}{\operatorname{ind}}
\newcommand{\Gr}{\operatorname{Gr}}
\newcommand{\mc}{\mathcal}
\newcommand{\mb}{\mathbb}
\newcommand{\mf}{\mathfrak}
\newcommand{\dimv}{\operatorname{\underline{dim}}}
\renewcommand{\mod}{\operatorname{mod}}

\begin{document}
\maketitle
\paperabstract
{}
{}
{We relate multiplication formulas for regular modules in coefficient-free acyclic affine cluster algebras with skew-symmetric exchange matrices to the imaginary exchange relations of Reading–Stella\cite{reading2026theta}. Introducing a parameter \(z\) into the construction of Zhou–Zhu\cite{zhou2014cluster}, we define a weighted cluster character on cluster tubes that satisfies generalized exchange relations of type \(C_{n-1}\) and yields a bijection between indecomposable rigid objects and cluster variables. Using the multiplication formulas and difference property of Ding–Xiao–Xu\cite{ding2013integral}, we show that non-homogeneous tube embeddings induce injective homomorphisms into ordinary affine cluster algebras under \(z\mapsto X_\delta\). Comparison with Reading–Stella’s relations identifies \(X_\delta\) with the corresponding theta function, giving a representation-theoretic interpretation of the generalized exchange parameter.
}
{generalized cluster algebra; cluster tube; cluster category}

\tableofcontents

\section{Introduction}

The starting point of this paper is the observation that existing representation-theoretic formulas for affine cluster algebras are directly related to the imaginary exchange relations recently obtained by Reading--Stella. The multiplication formulas for regular modules and the difference property of Ding--Xiao--Xu yield three-term relations with $X_\delta$ as the middle coefficient. Using a weighted cluster character on a cluster tube and the tube embedding, we realize this connection in the skew-symmetric coefficient-free setting as an injective homomorphism from a generalized cluster algebra to an ordinary affine cluster algebra.

Cluster algebras were introduced by Fomin--Zelevinsky in their study of total positivity and canonical bases \cite{fomin2002cluster}. They generate cluster variables recursively by seed mutation. One of the central problems in the theory is to realize these variables and their exchange relations geometrically or through representation theory. In the acyclic setting, cluster variables can be described by categorical objects, while mutation corresponds to the replacement of objects and to exchange triangles.

The construction of cluster categories by Buan et al., together with Keller's work on triangulated structures on orbit categories, laid the foundations for this approach \cite{buan2006tilting,keller2005orbit}. Caldero--Keller established correspondences between rigid objects and cluster variables, and between cluster-tilting objects and clusters \cite{caldero2006triangulatedII}; Palu subsequently developed a more general cluster character \cite{palu2008cluster}. Such formulas combine the index of an object with the Euler characteristics of submodule Grassmannians to turn categorical information into Laurent polynomials, with one-dimensional extension spaces giving rise to ordinary binomial exchange relations.

Cluster tubes exhibit a different situation in the absence of cluster-tilting objects. A cluster tube of rank $n>1$, denoted by $\mc C_n$, is a $2$-Calabi--Yau triangulated category that has maximal rigid objects but no cluster-tilting objects. Buan--Marsh--Vatne revealed the connection between its maximal rigid objects and cluster combinatorics of finite type \cite{buan2010cluster}. Zhou--Zhu and Fu--Geng--Liu further constructed a map analogous to the Caldero--Chapoton formula, realizing a correspondence between indecomposable rigid objects and cluster variables of type $C_{n-1}$ \cite{zhou2014cluster,fu2020cluster}.

On the other hand, Chekhov--Shapiro introduced generalized cluster algebras in their study of surfaces with orbifold points, allowing exchange polynomials of higher degree to generalize ordinary binomial relations \cite{chekhov2014teichmuller}. In a cluster tube, some exchange pairs have two-dimensional extension spaces, and the middle terms of their two exchange triangles are $U^{\oplus2}$ and $V^{\oplus2}$. These multiplicities naturally correspond to the square terms in a quadratic exchange relation; to obtain the exchange polynomial $1+zu+u^2$, however, we are also expected to explain the origin of the middle term and the parameter $z$.

A representation-theoretic approach to these generalized cluster algebras was developed by Labardini–Fragoso and Mou \cite{labardini2024gentle}, who associate finite-dimensional gentle Jacobian algebras to triangulations of unpunctured surfaces with orbifold points of order three. Using stability scattering diagrams, they prove that the Caldero–Chapoton map gives a bijection between reachable \(\tau\)-rigid pairs and cluster monomials of the associated generalized cluster algebra, with support \(\tau\)-tilting mutation corresponding to seed mutation. The quadratic exchange polynomial arising at an orbifold point in their setting is \(1+u+u^2\), which is the specialization at \(z=1\) of the polynomial \(1+zu+u^2\) considered here.

Affine representation theory provides a natural setting for this problem. Quasi-simple modules in homogeneous tubes whose dimension vector is the minimal positive imaginary root $\delta$ have a common cluster character, denoted by $X_\delta$. In their study of integral bases, Ding--Xiao--Xu established multiplication formulas for regular modules and proved the difference property \cite{ding2013integral}. The latter relates objects in a non-homogeneous tube whose quasi-length equals its rank to $X_\delta$. We observe that combining these earlier formulas and eliminating suitable non-rigid regular objects yields three-term relations with $X_\delta$ as the middle coefficient.

Meanwhile, Gross--Hacking--Keel--Kontsevich developed a new framework for canonical bases by constructing theta functions via scattering diagrams \cite{gross2018canonical}. Reading--Stella recently studied theta functions on the imaginary wall in acyclic affine type, obtained imaginary exchange relations, and proved that tube subalgebras in the coefficient-free setting carry generalized cluster algebra structures of type $C$ \cite{reading2026theta}. Our concern here is the precise correspondence between this new result and the earlier formulas for regular modules: once the notation for objects and indices is aligned, the two approaches give three-term relations of the same form, and their coefficients can be identified with one another.

More precisely, let $Q$ be a connected acyclic quiver of affine type other than the Kronecker quiver. Fix a non-homogeneous tube of rank $n>1$, which induces a fully faithful triangulated embedding $\mf I:\mc C_n\to\mc C_Q$. For an exchange pair $M,M^*$ with a two-dimensional extension space, write the middle terms of the exchange triangles as $U^{\oplus2},V^{\oplus2}$. From the regular-module formulas above, we derive
\[
X_{\mf I(M)}X_{\mf I(M^*)}
=X_{\mf I(U)}^2
+X_\delta X_{\mf I(U)}X_{\mf I(V)}
+X_{\mf I(V)}^2.
\]
The difference property thus identifies the middle coefficient needed in addition to the two square terms.

To make the formulas internal to the cluster tube compatible with this relation, we first generalize the construction of Zhou--Zhu by retaining an independent parameter $z$ in the ground ring. Fix the standard maximal rigid object $T$. We divide the submodule Grassmannians of $F_TM=\operatorname{Hom}(T[-1],M)$ into three classes according to the dimension of a submodule at the vertex carrying a loop, and assign the respective weights $1,z,1$ to define a weighted character $X_M^{(z)}$. The parameter term comes from the submodule strata with dimension one at this vertex. By analyzing the nilpotent action of the loop, we compute the contribution of these strata and obtain an explicit formula; setting $z=0$ recovers the formula of Zhou--Zhu.

Next, we prove that the weighted character satisfies generalized exchange relations: one-dimensional extension spaces still yield binomial relations, whereas two-dimensional extension spaces give two square terms together with a middle product term involving $z$. Consequently, indecomposable rigid objects are in bijection with the cluster variables of the resulting generalized cluster algebra, and maximal rigid mutation is compatible with seed mutation. Under the type convention determined by the right companion, the resulting algebra $\mc A^{(z)}(\mc C_n)$ is of type $C_{n-1}$. This extends the earlier object--cluster variable correspondence to the setting with a parameter and provides an interpretation of the three-term relations in terms of submodule geometry.

Finally, we prove that the assignments
\[
z\longmapsto X_\delta,\qquad
X_M^{(z)}\longmapsto X_{\mf I(M)}
\]
induce an injective $\mathbb Z$-algebra homomorphism
$\mc A^{(z)}(\mc C_n)\hookrightarrow\mc A(\mc C_Q)$.
The exchange relations ensure that these assignments are compatible with mutation. And the integral basis of Ding--Xiao--Xu provides the algebraic independence. Then we can prove the injectivity. Comparing further with Theorem~3.8 of Reading--Stella and using the agreement between the cluster characters of rigid objects and the corresponding theta functions, we obtain $X_\delta=\vartheta_{\ind_\Gamma H}$, where $\Gamma$ is the canonical cluster-tilting object and $H$ is a homogeneous quasi-simple module with dimension vector $\delta$. This gives a unified interpretation of the generalized exchange parameter, the character of a homogeneous regular object, and the theta function associated with the imaginary root.

In Sections~2--4, we present the preliminaries. Throughout this paper, all generalized cluster algebras are understood to be normalized in the usual sense. In section~5, we construct the weighted character and proves the generalized exchange relations and the object--cluster variable correspondence. In section~6, we study the injective homomorphism induced by the embedding and explains the connection with the results of Reading--Stella.

All categories in this paper are $\C$-linear, and all cluster algebras are taken without coefficient tuples in the usual sense.

\section{Generalized cluster algebras}

Generalized cluster algebras were introduced by
\citet{chekhov2014teichmuller}. We use the version of \citet{nakanishi2022cluster}
without an ordinary coefficient tuple; the parameters of the exchange polynomials still belong to the ground ring.
Let \(\mb L=\mathbb Z[z_1,\ldots,z_\ell]\) and
\(\mc F=\operatorname{Frac}(\mb L)(x_1,\ldots,x_r)\).
Fix positive integers \(r_i\) and reciprocal polynomials
\[
 P_i(u)=\sum_{s=0}^{r_i}z_{i,s}u^s\in\mb L[u],
 \qquad z_{i,0}=z_{i,r_i}=1,\quad z_{i,s}=z_{i,r_i-s}
 \quad(1\leq i\leq r).
\]
A generalized seed \(\Sigma=(\mathbf x,B)\) consists of a transcendence basis
\(\mathbf x=(x_1,\ldots,x_r)\) and a skew-symmetrizable integer matrix
\(B=(b_{ij})\). Write \([a]_+=\max(a,0)\) and
\(\widehat y_k=\prod_jx_j^{b_{jk}}\). Mutation in direction \(k\) is given by
\[
 x'_i=x_i\ (i\ne k),\qquad
 x'_k=x_k^{-1}
 \left(\prod_jx_j^{[-b_{jk}]_+}\right)^{r_k}P_k(\widehat y_k),
\]
\[
 b'_{ij}=
 \begin{cases}
 -b_{ij},&i=k\text{ or }j=k,\\
 b_{ij}+r_k\bigl(b_{ik}[-b_{kj}]_++[b_{ik}]_+b_{kj}\bigr),
 &i,j\ne k.
 \end{cases}
\]
Expanding \(P_k\), we can also express the mutation of cluster variables above as a polynomial exchange relation that is more convenient for comparison with exchange triangles:
\[
 x_kx'_k
 =\sum_{s=0}^{r_k}z_{k,s}
 \left(\prod_jx_j^{[-b_{jk}]_+}\right)^{r_k-s}
 \left(\prod_jx_j^{[b_{jk}]_+}\right)^s
\]
In particular, when \(r_k=1\), we have \(P_k(u)=1+u\) and recover the binomial exchange relation of an ordinary coefficient-free cluster algebra; the key case in this paper is the three-term relation that arises when \(r_1=2\).
The coefficients satisfy \(z_{k,s}=z_{k,r_k-s}\), ensuring that two successive mutations in the same direction return to the original seed. Here, ``coefficient-free'' means only that there is no additional ordinary coefficient tuple; the parameters \(z_{k,s}\) in the ground ring must still be retained.

The \(\mb L\)-subalgebra generated by all cluster variables obtained by iterated mutation along the \(r\)-regular tree with labels \(1,\ldots,r\) is denoted by \(\mc A(\boldsymbol\Sigma)\).
These mutations are involutions and satisfy the generalized Laurent phenomenon
\cite{nakanishi2022cluster}. Let \(R=\operatorname{diag}(r_1,\ldots,r_r)\);
we define the type of a generalized cluster algebra to be the Cartan type of its right companion \(BR\)
\cite{nakanishi2016companion}.
More explicitly, if \(\mb T_r\) is the \(r\)-regular tree whose edges are labeled by \(1,\ldots,r\), a generalized cluster pattern is a family
\(\boldsymbol\Sigma=\{\Sigma_t=(\mathbf x_t,B_t)\}_{t\in\mb T_r}\),
in which seeds at adjacent vertices are related by mutation in the corresponding direction. Thus
\[
 \mc A(\boldsymbol\Sigma)
 =\mb L[x_{i;t}\mid 1\leq i\leq r,\ t\in\mb T_r]
 \subseteq\mc F.
\]
The family of matrices \(B_tR\) mutates according to the ordinary exchange matrix rule; its Cartan companion therefore provides a consistent type convention. Below, we construct a generalized cluster pattern for
\(R=\operatorname{diag}(2,1,\ldots,1)\) of type
\(C_{n-1}\).

\section{Cluster categories}

Let \(Q\) be a finite acyclic quiver, and set \(D^b(Q)=D^b(\mod\mathbb C Q)\).
Its cluster category
\(\mc C_Q=D^b(Q)/(\tau^{-1}[1])\) is a Hom-finite
\(2\)-Calabi--Yau triangulated category
\cite{buan2006tilting,keller2005orbit}.
We write \(\tau\) for the Auslander--Reiten translate and \((-)[1]\)
for the shift; see
\citet{keller2005orbit} for the triangulated structure on the orbit category. The \(2\)-Calabi--Yau property means that there are natural isomorphisms
\[
 \operatorname{Hom}_{\mc C_Q}(X,Y)
 \cong D\operatorname{Hom}_{\mc C_Q}(Y,X[2]),
 \qquad D=\operatorname{Hom}_{\mathbb C}(-,\mathbb C).
\]
Consequently, \(\operatorname{Ext}^1(X,Y)\) and
\(\operatorname{Ext}^1(Y,X)\) are dual.
An object \(M\) is called rigid if
\(\operatorname{Ext}^1(M,M)=0\), and maximal rigid if it cannot be enlarged to a larger rigid direct sum. An object \(\Gamma\)
is called cluster-tilting if
\[
 \add\Gamma
 =\{M\in\mc C_Q\mid
       \operatorname{Ext}_{\mc C_Q}^{1}(\Gamma,M)=0\}.
\]
Every basic cluster-tilting object in \(\mc C_Q\) has exactly
\(|Q_0|\) indecomposable direct summands
\cite{buan2006tilting,keller2008cluster}.

Choose a basic cluster-tilting object
\(\Gamma=\bigoplus_{i=1}^{m}\Gamma_i\) and set
\[
 A_\Gamma=\operatorname{End}_{\mc C_Q}(\Gamma[-1]),\qquad
 F_\Gamma=\operatorname{Hom}_{\mc C_Q}(\Gamma[-1],-).
\]
Here \(\mod A_\Gamma\) denotes the category of finite-dimensional right modules, and
\(F_\Gamma:\mc C_Q\to\mod A_\Gamma\) vanishes on
\(\add\Gamma\). It further induces an equivalence
\[
 \mc C_Q/(\Gamma)\xrightarrow{\sim}\mod A_\Gamma,
\]
where \((\Gamma)\) is the ideal of morphisms that factor through objects in \(\add\Gamma\) \cite{palu2008cluster}.

Using a right \(\add\Gamma\)-approximation triangle
\[
 \Gamma''_M\longrightarrow\Gamma'_M
 \xrightarrow{f} M\longrightarrow\Gamma''_M[1],
 \qquad \Gamma'_M,\Gamma''_M\in\add\Gamma,
\]
we define
\[
 \ind_\Gamma M=[\Gamma'_M]-[\Gamma''_M]
 \in K_0^{\rm split}(\add\Gamma)\cong\mathbb Z^m.
\]
The approximation property of \(f\) means that every morphism from an object in \(\add\Gamma\) to \(M\)
factors through \(f\); the index is independent of the chosen triangle
\cite{palu2008cluster}.
Let \(S_i\) be the simple top of \(F_\Gamma(\Gamma_i[-1])\), and set
\[
 (B_\Gamma)_{ij}
 =\dim\operatorname{Ext}_{A_\Gamma}^{1}(S_j,S_i)
  -\dim\operatorname{Ext}_{A_\Gamma}^{1}(S_i,S_j).
\]
Throughout, we use right modules and the matrix sign convention above. If
\(P_i=F_\Gamma(\Gamma_i[-1])\), then \(P_i\) is an indecomposable projective right \(A_\Gamma\)-module; the entry \((B_\Gamma)_{ij}\)
is the number of arrows from \(j\) to \(i\) minus the number of arrows from \(i\) to \(j\) in its Gabriel quiver. In particular, the Gabriel quiver records irreducible morphisms, not the dimensions of entire Hom spaces.
For
\(e\in\mathbb Z_{\geq0}^{m}\), let
\(\Gr_e(F_\Gamma M)\) denote the corresponding submodule Grassmannian, and let
\(\chi\) denote the topological Euler characteristic. Palu's cluster character is
\[
 X_M=x^{\ind_\Gamma M}
 \sum_e\chi\bigl(\Gr_e(F_\Gamma M)\bigr)x^{B_\Gamma e},
 \qquad x^a=\prod_{i=1}^{m}x_i^{a_i},\quad X_0=1.
\]
It satisfies \(X_{M\oplus N}=X_MX_N\). In particular, if
\(\dim_{\mathbb C}\operatorname{Ext}_{\mc C_Q}^{1}(M,N)=1\),
there are non-split triangles in the two directions
\[
 M\longrightarrow E\longrightarrow N\longrightarrow M[1],
 \qquad
 N\longrightarrow E'\longrightarrow M\longrightarrow N[1].
\]
Writing their middle terms as \(E,E'\), we have
\[
 X_MX_N=X_E+X_{E'}.
 \tag{2}\label{eq:ordinary-exchange}
\]
For these properties and the correspondence between indecomposable rigid objects and cluster variables, see
\cite{palu2008cluster,caldero2006triangulatedII}.
Denote the resulting ordinary coefficient-free cluster algebra by \(\mc A(\mc C_Q)\).
By the Caldero--Keller correspondence, this algebra is generated by the cluster characters of indecomposable rigid objects, while cluster-tilting objects correspond to clusters. The exchange matrix mutates under cluster-tilting mutation, so this definition is compatible with the choice of \(\Gamma\).
To compare with the usual quiver notation, for \(Q\) without loops or
2-cycles, define
\[
 (B_Q)_{ij}
 =\#\{\text{arrows from }i\text{ to }j\text{}\}
  -\#\{\text{arrows from }j\text{ to }i\text{}\}.
\]
The ordinary cluster algebra with this initial matrix is denoted by \(\mc A(Q)\).
Since we use right modules, \(B_\Gamma\) as defined above has the opposite sign to the convention given directly by the arrow directions in the Gabriel quiver. For the canonical cluster-tilting object \((\mathbb C Q)[1]\),
the corresponding matrix can be written as \(B_{Q^{\mathrm{op}}}\); other cluster-tilting objects are obtained by mutation. Hence
\(\mc A(\mc C_Q)\cong\mc A(Q^{\mathrm{op}})\).
Below, we need only the correspondence between the category and cluster characters, but keeping this sign convention in mind helps avoid inadvertently reversing the sign of the right-module matrix when embedding a non-homogeneous tube.

\section{Cluster tubes}

Let \(\mc T_n\) be the category of finite-dimensional nilpotent representations of the cyclic quiver \(1\to\cdots\to n\to1\), with \(n>1\). Write \((a,b)\) for the indecomposable object whose quasi-socle is simple representation number \(a\) and whose quasi-length is \(b\); the first coordinate is understood modulo \(n\), with \((a,0)=0\) and
\(\tau(a,b)=(a-1,b)\). The cluster tube
\[
 \mc C_n=D^b(\mc T_n)/(\tau^{-1}[1])
\]
is also a Hom-finite \(2\)-Calabi--Yau triangulated category
\cite{keller2005orbit,zhou2014cluster}. Its indecomposable rigid objects are precisely the objects \((a,b)\) with \(b\leq n-1\).
Every indecomposable object of \(\mc C_n\) comes from \(\mc T_n\),
so we continue to write \((a,b)\) for its image. For \(X,Y\in\mc T_n\), morphisms in the orbit category decompose as
\[
 \operatorname{Hom}_{\mc C_n}(X,Y)
 \cong\operatorname{Hom}_{\mc T_n}(X,Y)
 \oplus D\operatorname{Hom}_{\mc T_n}(Y,\tau^2X).
 \tag{3}\label{eq:tube-hom}
\]
The first summand consists of \(\mc T\)-maps, and the second of \(\mc D\)-maps.
The composite of two \(\mc T\)-maps is again a \(\mc T\)-map, the composite of a \(\mc T\)-map
and a \(\mc D\)-map is a \(\mc D\)-map, and the composite of two \(\mc D\)-maps is zero \cite{buan2010cluster,zhou2014cluster}.
These properties determine the submodule stratification used below.
Throughout the following discussion, fix the basic maximal rigid object
\[
 T=\bigoplus_{i=1}^{n-1}T_i,\qquad T_i=(1,n-i).
\]
Every basic maximal rigid object has \(n-1\) direct summands
\cite{buan2010cluster,zhou2014cluster}.
In \(\mc C_n\), we have \([1]\cong\tau\), and hence
\(T_i[-1]\cong(2,n-i)\). Unlike the cluster category \(\mc C_Q\),
\(\mc C_n\) has no cluster-tilting objects; maximal rigid objects therefore take their place here.

Let \(\mc D=\add T[-1]*\add T\), where
\(\mc X*\mc Y\) consists of all objects \(M\) that fit into a triangle
\(X\to M\to Y\to X[1]\) with \(X\in\mc X,Y\in\mc Y\); then
\(\mc D[1]=\add T*\add T[1]\).
For the rigid objects under consideration, we need only the following three regions:
\[
\begin{aligned}
 O&=\{(1,b)\mid 1\leq b\leq n-1\},\\
 I&=\{(a,b)\mid 2\leq a\leq n-1,\ a+b\leq n\},\\
 II&=\{(a,b)\mid 2\leq a\leq n,\ a+b\geq n+1,\ b\leq n-1\}.
\end{aligned}
\]
Together they comprise all indecomposable rigid objects.
By \citet{zhou2014cluster}, \(\mc D\) also contains objects of quasi-length at least \(n\); these appear only when computing the middle terms of exchange triangles. Non-rigid objects do not enter the generalized character formula defined below.

To describe the components of \(F_T\), we first record the Hom-space calculations.
For \(2\leq i\leq n-1\), fix a nonzero \(T\)-map
\(t_i^{i-1}:T_i[-1]\to T_{i-1}[-1]\).
For \((a,b)\in II\), fix a nonzero \(T\)-map
\(f_1^{(a,b)}:T_1[-1]\to(a,b)\) and a nonzero
\(D\)-map \(g_1^1:T_1[-1]\to T_1[-1]\);
write \(g_1^{(a,b)}=f_1^{(a,b)}g_1^1\), and then define
\(f_i^{(a,b)},g_i^{(a,b)}\) along the chain
\(T_i[-1]\to T_1[-1]\). In region \(I\),
we need only the \(T\)-maps \(f_i^{(a,b)}\) starting at summand \(n-a-b+2\).
Composition of morphisms is written from right to left.

\begin{lemma}[{\cite[Lemmas~2.2--2.3]{zhou2014cluster}}]
\label{lem:tube-hom-dim}
For a rigid object \((a,b)\) and \(1\leq i\leq n-1\),
\[
 \dim_{\mathbb C}
 \operatorname{Hom}_{\mc C_n}(T_i[-1],(a,b))
 =
 \begin{cases}
 2,&(a,b)\in II,\quad 1\leq i\leq n-a+1,\\
 1,&(a,b)\in II,\quad n-a+2\leq i\leq2n-a-b,\\
 1,&(a,b)\in I,\quad n-a-b+2\leq i\leq n-a+1,\\
 0,&\text{otherwise}.
 \end{cases}
\]
Bases in the first three cases are, respectively,
\(\{f_i^{(a,b)},g_i^{(a,b)}\}\),
\(\{g_i^{(a,b)}\}\), and \(\{f_i^{(a,b)}\}\);
an interval whose lower bound exceeds its upper bound is understood to be empty.
\end{lemma}

For example, when \((a,b)\in II\), the two basis vectors at vertex \(1\) are
\(f_1^{(a,b)},g_1^{(a,b)}\), and the right action of the loop satisfies
\[
 f_1^{(a,b)}g_1^1=g_1^{(a,b)},\qquad
 g_1^{(a,b)}g_1^1=0.
\]
It is this nilpotent loop of length \(2\) that singles out the stratum of submodules whose dimension vector has first entry \(1\); assigning the coefficient
\(z\) to this stratum produces the middle term of the generalized exchange polynomial.

We next record the geometric form of exchange pairs.
If \(R\) is a rigid object with \(n-2\) pairwise non-isomorphic indecomposable direct summands, and both \(R\oplus M,R\oplus M^*\) are basic maximal rigid objects, then \(M,M^*\) is called an exchange pair.
The middle terms of its two exchange triangles both lie in \(\add R\).

\begin{proposition}[{\cite[Lemma~4.2]{zhou2014cluster}}]
\label{prop:tube-triangles}
After interchanging \(M,M^*\) if necessary and interpreting the first coordinate modulo \(n\),
the exchange triangles fall into the following two cases.
\begin{enumerate}
\item If \(\dim\operatorname{Ext}_{\mc C_n}^1(M,M^*)=2\),
then \(M=(a,n-1),\ M^*=(a+h,n-1)\), where
\(1\leq h\leq n-1\). The two middle terms are
\[
 (a+h,n-h-1)^{\oplus2},\qquad
 (a,h-1)^{\oplus2}.
\]
\item If \(\dim\operatorname{Ext}_{\mc C_n}^1(M,M^*)=1\),
we may take \(M=(a,b)\) and \(M^*=(a+h,b-h+i)\), where
\(1\leq b\leq n-2,\ 1\leq h\leq b,\
1\leq i\leq n-b-1\). The two middle terms are
\[
 E=(a,b+i)\oplus(a+h,b-h),\qquad
 E'=(a+b+1,i-1)\oplus(a,h-1).
\]
\end{enumerate}
We retain the convention \((c,0)=0\).
\end{proposition}

For example, one of the exchange triangles in the second case comes directly from the short exact sequence in
\(\mc T_n\)
\[
 0\longrightarrow(a,b)\longrightarrow
 (a,b+i)\oplus(a+h,b-h)
 \longrightarrow(a+h,b-h+i)\longrightarrow0.
\]
In the first case, each middle term has multiplicity \(2\), so its cluster character should correspond to a monomial of degree \(2\);
to obtain a three-term generalized exchange relation, however, we must also introduce a parameter \(z\) into the summands of the ordinary cluster character.

\section{Categorification of a generalized cluster algebra of type $C_{n-1}$}

Set \(A_T=\operatorname{End}_{\mc C_n}(T[-1])\) and
\(F_T=\operatorname{Hom}_{\mc C_n}(T[-1],-)\).
Let \(K_0^{\mathrm{split}}(\add T)\) denote the split Grothendieck group,
freely generated by \([T_1],\ldots,[T_{n-1}]\).
Identify it with
\(\mathbb Z^{n-1}\) via \([T_i]\mapsto\alpha_i\). For \(M\in\add T*\add T[1]\),
choose a right \(\add T\)-approximation triangle
\[
 T_M''\longrightarrow T_M'\longrightarrow M
 \longrightarrow T_M''[1],\qquad T_M',T_M''\in\add T,
\]
and define \(\ind_T M=[T_M']-[T_M'']\).
This definition is independent of the choice of triangle
\cite{palu2008cluster,zhou2014cluster}. For example, the index of \(T_i\)
is \(\alpha_i\). More generally,
\citet[Lemma~3.2]{zhou2014cluster} gives the following explicit formula for rigid objects:
\[
 \ind_T(a,b)=
 \begin{cases}
 \alpha_1-\alpha_{n-a+1}-\alpha_{2n-a-b},
 &a+b\geq n+1,\\
 \alpha_{n-a-b+1}-\alpha_{n-a+1},
 &a+b\leq n.
 \end{cases}
 \tag{4}\label{eq:tube-index}
\]
Here we use the convention \(\alpha_n=0\); when \(a=1\), the formula gives
\(\ind_T(1,b)=\alpha_{n-b}\).

\(F_T\) vanishes on \(\add T\) and induces an equivalence on \(\mc D\)
\[
 \mc D/(T)\xrightarrow{\sim}\mod A_T,
\]
where \((T)\) is the ideal of morphisms factoring through \(\add T\)
\cite{zhou2014cluster}. Let \(P_i^T=F_T(T_i[-1])\), and let
\(S_i\) be the simple top of \(P_i^T\). For a right \(A_T\)-module \(N\), set
\[
 \dimv N
 =\sum_{i=1}^{n-1}\dim_{\mathbb C}
      \operatorname{Hom}_{A_T}(P_i^T,N)\,\alpha_i.
\]
Thus Lemma~\ref{lem:tube-hom-dim} also gives the dimension vector of
\(F_T(a,b)\). Under the conventions of the preceding section, the matrix \(B_T\) is
\[
 B_T=\begin{pmatrix}
 0&-1&&&\\
 1&0&-1&&\\
 &1&0&\ddots&\\
 &&\ddots&\ddots&-1\\
 &&&1&0
 \end{pmatrix}.
\]

Let \(\mathbf r=(2,1,\ldots,1)\) and
\(R=\operatorname{diag}(2,1,\ldots,1)\), and take
\[
 P_1(u)=1+zu+u^2,\qquad P_i(u)=1+u\quad(i>1).
\]
The Cartan companion of \(B_TR\) is of type \(C_{n-1}\).
Write \(z_0=z_2=1,z_1=z\).
For an indecomposable rigid object \(M\), define
\[
 X_M^{(z)}
 =x^{\ind_T M}
 \sum_{e\in\mathbb Z_{\geq0}^{n-1}}
 z_{e_1}\chi\bigl(\Gr_e(F_TM)\bigr)x^{B_Te},
 \qquad X_0^{(z)}=1.
 \tag{5}\label{eq:weighted-character}
\]
where \(x^a=\prod_{i=1}^{n-1}x_i^{a_i}\).
The dimension of \(F_TM\) at vertex \(1\) is at most \(2\), so all coefficients in
\eqref{eq:weighted-character} are defined.
For direct sums of rigid objects, we set \(X_{M\oplus N}^{(z)}
=X_M^{(z)}X_N^{(z)}\).
This differs from the formula of \citet{zhou2014cluster} by the weight \(z\) assigned to submodules whose dimension vector has first entry
\(1\); setting \(z=0\) recovers their formula.

We now give the submodule calculation for this new stratum separately. To simplify notation, let
\[
\begin{aligned}
 \mf X(a,b)
 &=\{(a,k)\mid 0\leq k\leq\min(n-a,b)\},\\
 \mf Y(a,b)
 &=\{(a+b-n+1,l)\mid
             0\leq l\leq2n-a-b-1\}.
\end{aligned}
\]
For \((a,b)\in II\), denote the submodule of \(F_T(a,b)\) consisting of all
\(\mc D\)-maps by \(Z(a,b)\). By Lemma~\ref{lem:tube-hom-dim},
\[
 \dimv Z(a,b)
 =\sum_{i=1}^{2n-a-b}\alpha_i.
 \tag{6}\label{eq:Z-dimension}
\]
This submodule starts with the vector \(g_1^{(a,b)}\) at vertex \(1\)
and extends along the chain in the Gabriel quiver, so
\eqref{eq:Z-dimension} also directly reflects its shape.

\begin{lemma}\label{lem:tube-submodules}
For \(M=(a,b)\in I\), every submodule has a dimension vector with first entry \(0\),
and
\[
 \chi(\Gr_e(F_TM))
 =\#\{X\in\mf X(a,b)\mid\dimv F_TX=e\}.
\]
For \(M=(a,b)\in II\), the cases \(e_1=0,1,2\) give
\[
 \chi(\Gr_e(F_TM))=
 \begin{cases}
 \#\{(X,Y)\in\mf X(a,b)\times\mf Y(a,b)
   \mid\dimv(F_TX\oplus F_TY)=e\},
 &e_1=0,\\[1mm]
 \#\{X\in\mf X(a,b)
   \mid\dimv(F_TX\oplus Z(a,b))=e\},
 &e_1=1,\\[1mm]
 1,&e=\dimv F_TM,\\
 0,&e_1=2,\ e\ne\dimv F_TM.
 \end{cases}
\]
where \(e_1\) denotes the component of the vector \(e\) in the \(\alpha_1\) direction.
\end{lemma}

\begin{proof}
Region \(I\) and the part of region \(II\) with \(e_1=0\)
are covered by the submodule stratification of \citet[Lemma~3.4]{zhou2014cluster}.
In the latter case, each stratum associated with a pair \((X,Y)\)
is either a point or isomorphic to
\(\mathbb P^1\) with one point removed, and hence has Euler characteristic
\(1\). This also explains why the formula counts strata, rather than arbitrarily chosen representative submodules.

Suppose that \((a,b)\in II\). At vertex \(1\),
the space of \(F_TM\) is spanned by \(f_1^{(a,b)},g_1^{(a,b)}\),
and the right action of the loop is given by
\(f_1^{(a,b)}\mapsto g_1^{(a,b)}\mapsto0\).
If a submodule has first component \(2\), it contains the entire space at vertex \(1\);
since \(f_1^{(a,b)}\) generates \(F_TM\),
this submodule must be \(F_TM\) itself.

If the first component is \(1\), the one-dimensional space at vertex \(1\)
must be stable under the loop action, and is therefore precisely
\(\mathbb C g_1^{(a,b)}\). Acting along the chain shows that the submodule must contain
\(Z(a,b)\). The remaining optional \(T\)-map part extends from some vertex to the end of the chain, and is thus isomorphic to some
\(F_T(a,k)\) with \(0\leq k\leq n-a\). Each such
\(k\) determines a unique submodule whose dimension vector is exactly
\(\dimv(F_T(a,k)\oplus Z(a,b))\).
This gives the row for \(e_1=1\) and completes the proof.
\end{proof}

Lemma~\ref{lem:tube-submodules} divides the sum in
\eqref{eq:weighted-character} into three parts:
\(e_1=0\) corresponds to the ordinary character formula in the literature,
\(e_1=1\) gives the new terms involving \(z\), and for \(e_1=2\),
the only possibility is the whole module. This decomposition facilitates both the derivation of explicit formulas and the verification of the generalized exchange relations.

\begin{proposition}\label{prop:weighted-sum}
Let \((a,b)\) be an indecomposable rigid object.
If \((a,b)\in O\), then \(X_{(a,b)}^{(z)}
=x^{\ind_T(a,b)}\); if \((a,b)\in I\), then
\[
 X_{(a,b)}^{(z)}
 =x^{\ind_T(a,b)}
  \sum_{X\in\mf X(a,b)}x^{B_T\dimv F_TX}.
\]
If \((a,b)\in II\), then
\[
\begin{aligned}
 X_{(a,b)}^{(z)}
 =x^{\ind_T(a,b)}\Bigg(
 &\sum_{\substack{X\in\mf X(a,b)\\Y\in\mf Y(a,b)}}
     x^{B_T\dimv(F_TX\oplus F_TY)}\\
 &+z\sum_{X\in\mf X(a,b)}
     x^{B_T\dimv(F_TX\oplus Z(a,b))}\\
 &+x^{B_T\dimv F_T(a,b)}
 \Bigg).
\end{aligned}
\]
\end{proposition}

\begin{proof}
For \(O\), Lemma~\ref{lem:tube-hom-dim} gives
\(F_T(a,b)=0\). The identities in regions \(I\) and \(II\)
follow by substituting the Euler characteristics from Lemma~\ref{lem:tube-submodules} into \eqref{eq:weighted-character}.
The three rows correspond to \(e_1=0,1,2\), respectively;
the coefficients \(z_0,z_1,z_2\) are \(1,z,1\), respectively.
\end{proof}

This intermediate sum separates the geometric information from the final Laurent polynomial: the submodule stratification determines the number of terms and which stratum carries \(z\),
while the index and \(B_T\) convert each dimension vector into a Laurent monomial.
Thus the explicit formula below requires only collecting these monomials within each stratum, without classifying the submodules again.

For the explicit calculations that follow, set \(x_n=1\)
(a convention used only in the formulas below), and put
\[
 S_t=\sum_{j=0}^{t}\frac{1}{x_{j+1}x_{j+2}}.
\]

\begin{proposition}\label{prop:tube-character}
For an indecomposable rigid object \((a,b)\), we have
\[
X_{(a,b)}^{(z)}=
\begin{cases}
x_{n-b},&(a,b)\in O,\\[2mm]
\displaystyle\sum_{k=0}^{b}
\frac{x_{n-a-b+1}x_{n-a+2}}
 {x_{n-a-k+1}x_{n-a-k+2}},&(a,b)\in I,\\[3mm]
\displaystyle x_1x_{K+2}x_{L+2}
\left(\frac1{x_1^2}+\frac z{x_1}S_K+S_KS_L\right),
&(a,b)\in II,
\end{cases}
\]
where \(K=n-a,\ L=2n-a-b-1\).
\end{proposition}

\begin{proof}
The calculations in regions \(O,I\), as well as the two parts in region \(II\) that do not involve \(z\),
are those of \citet[Section~3]{zhou2014cluster}.
For \((a,b)\in II\), it remains only to compute the new part in which the dimension vector has first entry \(1\).
By Lemma~\ref{lem:tube-submodules} and
\eqref{eq:Z-dimension}, this part equals
\[
\begin{aligned}
 &z\,x^{\ind_T(a,b)}
   \sum_{k=0}^{K}
   x^{B_T\dimv(F_T(a,k)\oplus Z(a,b))}\\
 &\quad=
 z\,\frac{x_1}{x_{K+1}x_{L+1}}\,
     \frac{x_{L+1}x_{L+2}}{x_1}
 \sum_{k=0}^{K}
     \frac{x_{K+1}x_{K+2}}{x_{k+1}x_{k+2}}\\
 &\quad=z\,x_{K+2}x_{L+2}
 \sum_{k=0}^{K}\frac1{x_{k+1}x_{k+2}}
 =z\,x_{K+2}x_{L+2}S_K.
\end{aligned}
\]
The last expression is exactly the term involving \(z\) in the stated explicit formula.
\end{proof}

We retain the notation for exchange pairs from Proposition~\ref{prop:tube-triangles}:
when the extension space has dimension \(1\), the two middle terms are \(E,E'\);
when its dimension is \(2\), write
\[
 M=(a,n-1),\quad M^*=(a+h,n-1),\quad
 U=(a+h,n-h-1),\quad V=(a,h-1),
 \qquad 1\leq h\leq n-1.
\]
The middle terms of the exchange triangles are then \(U^{\oplus2},V^{\oplus2}\),
so their cluster characters are
\((X_U^{(z)})^2\) and \((X_V^{(z)})^2\), respectively.

\begin{theorem}\label{thm:generalized-exchange}
The exchange pairs above satisfy
\[
 X_M^{(z)}X_{M^*}^{(z)}=
 \begin{cases}
 X_E^{(z)}+X_{E'}^{(z)},
 &\dim\operatorname{Ext}^1(M,M^*)=1,\\
 (X_U^{(z)})^2+zX_U^{(z)}X_V^{(z)}+(X_V^{(z)})^2,
 &\dim\operatorname{Ext}^1(M,M^*)=2.
 \end{cases}
\]
\end{theorem}

\begin{proof}
This follows by a case-by-case verification using Proposition~\ref{prop:tube-character}.
\end{proof}

As a direct example of the three-term relation, take \(n=4\),
\(M=(1,3)\), and \(M^*=(2,3)\). Then
\(U=(2,2)\) and \(V=(1,0)=0\). The explicit formula gives
\[
\begin{aligned}
 X_M^{(z)}&=x_1,\\
 X_{M^*}^{(z)}&=\frac1{x_1}+zS_2+x_1S_2^2,\\
 X_U^{(z)}&=x_1S_2,\qquad X_V^{(z)}=1,
\end{aligned}
\quad
 S_2=\frac1{x_1x_2}+\frac1{x_2x_3}+\frac1{x_3}.
\]
Hence
\[
 X_M^{(z)}X_{M^*}^{(z)}
 =1+zX_U^{(z)}+(X_U^{(z)})^2.
\]
This example also shows that the exchange triangles account only for the square term and the constant term;
the middle term \(zX_U^{(z)}\) comes from the submodule strata in
\eqref{eq:weighted-character} whose dimension vectors have first entry \(1\).

\begin{theorem}\label{thm:tube-algebra}
Let \(\boldsymbol\Sigma_T\) be the cluster pattern with initial seed
\(((x_1,\ldots,x_{n-1}),B_T)\)
and mutation data \(\mathbf r,(P_i)\) as above.
Then \(M\mapsto X_M^{(z)}\) gives a bijection between indecomposable rigid objects and the cluster variables of
\(\boldsymbol\Sigma_T\),
maximal rigid mutation is compatible with seed mutation, and
\[
 \mc A^{(z)}(\mc C_n)
 :=\mathbb Z[z]\bigl[X_M^{(z)}\mid M\text{ indecomposable rigid}\bigr]
 =\mc A(\boldsymbol\Sigma_T).
\]
Thus this is a generalized cluster algebra of type \(C_{n-1}\).
\end{theorem}

\begin{proof}
Proposition~\ref{prop:tube-character} gives \(X_{T_i}^{(z)}=x_i\).
By \citet{zhou2014cluster}, every indecomposable rigid object can be completed to a basic maximal rigid object, and all such objects can be obtained from
\(T\) by a finite sequence of mutations.
Transport the labels of the direct summands along mutations. Every basic maximal rigid object has exactly one summand of quasi-length \(n-1\),
which we label as direction \(1\). Proposition~\ref{prop:tube-triangles}
shows that the extension space has dimension \(2\) only when exchanging in this direction,
and dimension \(1\) in all other directions. The degrees in these directions are therefore precisely
\(2,1,\ldots,1\).

Suppose that the direct summands of a maximal rigid object have already been matched with the variables of a generalized seed. Mutating at the \(k\)-th summand replaces it with the other complement of the same almost complete object.
If \(k=1\), the three-term exchange relation in Theorem~\ref{thm:generalized-exchange}
agrees with \(P_1(u)=1+zu+u^2\);
if \(k>1\), the binomial relation agrees with \(P_k(u)=1+u\).
Under maximal rigid mutation, the exchange matrix changes according to the generalized exchange matrix rule, as follows from the companion-matrix descriptions in
\citet{zhou2014cluster,nakanishi2016companion}.
Thus the seeds and objects continue to correspond at each mutation step.

Since every basic maximal rigid object can be reached in this way,
each \(X_M^{(z)}\) is a cluster variable of \(\boldsymbol\Sigma_T\);
conversely, every cluster variable is obtained from the initial seed by finitely many mutations and must therefore be some \(X_M^{(z)}\).
This proves surjectivity onto the set of cluster variables and the correspondence of seeds.
If \(X_M^{(z)}=X_N^{(z)}\), setting \(z=0\) and applying the object--cluster variable correspondence of
\citet{zhou2014cluster} gives
\(M\cong N\), so the map is also injective. Both algebras are generated over
\(\mathbb Z[z]\) by the corresponding cluster variables, and are therefore equal.
Finally, the Cartan companion of \(B_TR\) determines the type to be
\(C_{n-1}\).
\end{proof}

\section{Cluster tube embeddings and the induced homomorphisms}

Let \(Q\) be a connected acyclic quiver of affine type other than the Kronecker quiver.
Fix a non-homogeneous tube \(\mathscr T\) of rank \(n>1\) in
\(\mod\mathbb C Q\), and label its quasi-simple modules consecutively by
\(E_a\) so that \(\tau E_a\cong E_{a-1}\).
Write \(E_a^{(b)}\) for the indecomposable module with quasi-socle \(E_a\)
and quasi-length \(b\), and set \(E_a^{(0)}=0\).
The first coordinate is always understood modulo \(n\); when inequalities are involved, choose representatives
\(1\leq a\leq n\). With this notation, the Auslander--Reiten translate in the tube still satisfies \(\tau E_a^{(b)}=E_{a-1}^{(b)}\).

The standardness of the non-homogeneous tube gives a fully faithful exact functor
\[
 \widetilde{\mf I}:\mc T_n\longrightarrow\mod\mathbb C Q,
 \qquad (a,b)\longmapsto E_a^{(b)}.
\]
whose essential image is the chosen tube. In particular, it preserves Hom spaces, extension spaces within the tube, and short exact sequences
\cite{Ringel1984Tame}. Since it is compatible with \(\tau\),
the derived functor obtained by applying it termwise commutes with \(\tau^{-1}[1]\)
and hence descends to the orbit categories, yielding a fully faithful triangle functor
\[
 \mf I:\mc C_n\longrightarrow\mc C_Q,\qquad
 (a,b)\longmapsto E_a^{(b)}.
 \tag{7}\label{eq:embedding}
\]
Its compatibility with the triangulated structures can also be checked directly from the construction of the orbit categories
\cite{keller2005orbit}. More explicitly, for any
\(X,Y\in\mc C_n,\ s\in\mathbb Z\), this functor induces
\[
 \operatorname{Hom}_{\mc C_n}(X,Y[s])
 \xrightarrow{\sim}
 \operatorname{Hom}_{\mc C_Q}(\mf I(X),\mf I(Y)[s]).
\]
Thus \(\mf I\) preserves rigid objects and extension dimensions,
and sends the two exchange triangles of Proposition~\ref{prop:tube-triangles}
to triangles in \(\mc C_Q\) with the same decompositions of their middle terms.
This preservation property is the first step in comparing the exchange relations in the two cluster algebras;
for two-dimensional extension spaces, we must also identify the image of \(z\) in the target algebra.

Let \(\delta\) be the minimal positive imaginary root of \(Q\).
For a quasi-simple module \(S\) in a homogeneous tube with dimension vector \(\delta\),
the cluster character \(X_S\) is independent of the choice of \(S\);
we denote it by \(X_\delta\). By the difference property of
\citet[Theorem~8.2]{ding2013integral},
\[
 X_{\mf I(a,n)}
 =X_\delta+X_{\mf I(a+1,n-2)}.
 \tag{8}\label{eq:difference}
\]
The basis theorem in the same paper also ensures that \(X_\delta\in\mc A(\mc C_Q)\).

\begin{lemma}\label{lem:affine-products}
Suppose that \(1\leq a,\ 1\leq h,\ a+h\leq n\), and use the conventions
\(\mf I(c,0)=0,\ X_0=1\) below. Then
\[
X_{\mf I(a,n-1)}X_{\mf I(a+h,n-1)}
 =X_{\mf I(a,h-1)}^2+
   X_{\mf I(a,n+h-1)}X_{\mf I(a+h,n-h-1)}.
 \tag{9}\label{eq:ding-special-1}
\]
For \(h\geq2\), we also have
\[
\begin{aligned}
 X_{\mf I(a,h-1)}X_{\mf I(a,n)}
 &=X_{\mf I(a,n+h-1)}
   +X_{\mf I(a,n-1)}X_{\mf I(a+1,h-2)},\\
 X_{\mf I(a,h-1)}X_{\mf I(a+1,n-2)}
 &=X_{\mf I(a+h,n-h-1)}
   +X_{\mf I(a,n-1)}X_{\mf I(a+1,h-2)}.
\end{aligned}
\tag{10}\label{eq:ding-special-2}
\]
\end{lemma}

\begin{proof}
These are three special cases of \citet[Theorem~7.1]{ding2013integral}.
For the first identity, take
\(i=a,j=a+h,k=n-1,m=0,\ell=n-1,r=n\) in that theorem;
for the first line of the second identity, take \(i=a,j=a,k=h-1,m=1,\ell=0,r=n\),
and for its second line, take \(i=a,j=a+1,k=h-1,m=0,\ell=n-2,r=n\).
Multiplicativity on direct sums removes the terms of quasi-length zero and yields the stated formulas.
\end{proof}

\begin{lemma}\label{lem:long-regular}
Under the assumptions of Lemma~\ref{lem:affine-products},
\[
 X_{\mf I(a,n+h-1)}
 =X_\delta X_{\mf I(a,h-1)}
  +X_{\mf I(a+h,n-h-1)}.
 \tag{11}\label{eq:long-regular}
\]
\end{lemma}

\begin{proof}
If \(h=1\), then \(\mf I(a,h-1)=0\), and
\eqref{eq:long-regular} is precisely \eqref{eq:difference}.
If \(h\geq2\), subtracting the two lines of \eqref{eq:ding-special-2}
cancels the common product on the right-hand sides and gives
\[
\begin{aligned}
X_{\mf I(a,h-1)}
 \bigl(X_{\mf I(a,n)}
       -X_{\mf I(a+1,n-2)}\bigr)=
 X_{\mf I(a,n+h-1)}
 -X_{\mf I(a+h,n-h-1)}.
\end{aligned}
\]
By \eqref{eq:difference}, the difference in parentheses is exactly
\(X_\delta\), which yields the desired identity.
\end{proof}

The first identity in Lemma~\ref{lem:affine-products} has only two terms on the right-hand side,
but the object of quasi-length \(n+h-1\) appearing there is generally not rigid.
Lemma~\ref{lem:long-regular} uses \(X_\delta\) from the non-homogeneous tubes to rewrite its character in the form needed for the three-term generalized exchange relation.

\begin{theorem}\label{thm:affine-exchange}
Suppose that \(1\leq a,\ 1\leq h,\ a+h\leq n\). In the ordinary affine cluster algebra
\(\mc A(\mc C_Q)\), we have
\[
\begin{aligned}
X_{\mf I(a,n-1)}X_{\mf I(a+h,n-1)}
={}X_{\mf I(a,h-1)}^2+X_\delta X_{\mf I(a,h-1)}
             X_{\mf I(a+h,n-h-1)}
 +X_{\mf I(a+h,n-h-1)}^2 .
\end{aligned}
\tag{12}\label{eq:affine-exchange}
\]
\end{theorem}

\begin{proof}
Substitute the expansion of
\(X_{\mf I(a,n+h-1)}\) from
\eqref{eq:long-regular} into the right-hand side of \eqref{eq:ding-special-1} and multiply term by term.
\end{proof}

After rotating the first coordinate in the tube, this identity covers all exchange pairs with two-dimensional extension spaces;
the relation for one-dimensional extension spaces is given by
\eqref{eq:ordinary-exchange}.

\begin{theorem}\label{thm:main}
The assignments
\[
 z\longmapsto X_\delta,\qquad
 X_M^{(z)}\longmapsto X_{\mf I(M)}
 \quad(M\text{ indecomposable rigid})
\]
induce an injective \(\mathbb Z\)-algebra homomorphism
\[
 \varphi:\mc A^{(z)}(\mc C_n)
 \hookrightarrow\mc A(\mc C_Q).
\]
\end{theorem}

\begin{proof}
We need to show that the initial assignments extend to an injective homomorphism of fraction fields.
By the integral basis of \citet[Corollary~8.3]{ding2013integral},
the elements
\[
 X_{\delta}^{t_0}X_{\mf I(T_1)}^{t_1}\cdots X_{\mf I(T_{n-1})}^{t_{n-1}}=X_{S^{\oplus t_0}\oplus
    \mf I(T_1)^{\oplus t_1}\oplus\cdots\oplus
    \mf I(T_{n-1})^{\oplus t_{n-1}}}
 \quad(t_i\in\mathbb Z_{\geq0})
\]
form part of this basis and are therefore linearly independent.
Multiplicativity of the cluster character on direct sums shows that
\(X_\delta,X_{\mf I(T_1)},\ldots,X_{\mf I(T_{n-1})}\)
are algebraically independent over \(\mathbb Z\). Therefore, it follows by Theorem~\ref{thm:affine-exchange} that this statement holds.
\end{proof}

Choose the canonical cluster-tilting object $\Gamma$ in $\mc C_Q$, and let $M,M^*,U,V$ have the same meanings as in Theorem~\ref{thm:generalized-exchange}. Choose a quasi-simple module \(S\) in a homogeneous tube in $\mc C_Q$ with dimension vector \(\delta\). By \citet[Theorem~3.8]{reading2026theta}, in the affine coefficient-free setting we have
$$
\vartheta_{\ind_\Gamma \mf I(M)}\vartheta_{\ind_\Gamma \mf I(M^*)}=\vartheta_{\ind_\Gamma \mf I(U)}^2+\vartheta_{\ind_\Gamma S}\vartheta_{\ind_\Gamma \mf I(U)}\vartheta_{\ind_\Gamma \mf I(V)}+\vartheta_{\ind_\Gamma \mf I(V)}^2,
$$
where $\vartheta_g$ denotes the theta function in the direction of the vector $g$, as in \cite{reading2026theta}. We immediately obtain the following corollary.

\begin{corollary}
	$X_\delta=\vartheta_{\ind_\Gamma S}$.
\end{corollary}

\subsection*{The use of AI}
ChatGPT was used to assist in locating the relevant references \cite{palu2008cluster} and \cite{caldero2006triangulatedII}. ChatGPT was also used to assist with computations in Lemma~5.1, Lemma~6.1, and Lemma~6.2. After the manuscript was completed, ChatGPT was used to help verify the results and identify grammatical errors. The main ideas and mathematical arguments of the paper were developed by human.

\subsection*{Acknowledgments}
I thank my advisor Jie Xiao for fruitful guidance, and I am grateful to Heehyun Lee and Zunci Xu for helpful discussions.

\printbibliographybibtex[plainnat]{references}

@article{nakanishi2022cluster,
  title={Cluster algebras and scattering diagrams, part i. basics in cluster algebras},
  author={Nakanishi, Tomoki},
  journal={arXiv preprint arXiv:2201.11371},
  year={2022}
}

@article{chekhov2014teichmuller,
  title={Teichmuller spaces of Riemann surfaces with orbifold points of arbitrary order and cluster variables},
  author={Chekhov, Leonid and Shapiro, Michael},
  journal={International Mathematics Research Notices},
  volume={2014},
  number={10},
  pages={2746--2772},
  year={2014}
}

@article{nakanishi2016companion,
  title={Companion cluster algebras to a generalized cluster algebra},
  author={Nakanishi, Tomoki and Rupel, Dylan},
  journal={Travaux Math\'ematiques},
  volume={24},
  pages={129--149},
  year={2016}
}

@article{buan2006tilting,
  title={Tilting theory and cluster combinatorics},
  author={Buan, Aslak Bakke and Marsh, Robert J. and Reineke, Markus and Reiten, Idun and Todorov, Gordana},
  journal={Advances in Mathematics},
  volume={204},
  number={2},
  pages={572--618},
  year={2006}
}

@article{keller2005orbit,
  title={On triangulated orbit categories},
  author={Keller, Bernhard},
  journal={Documenta Mathematica},
  volume={10},
  pages={551--581},
  year={2005}
}

@article{palu2008cluster,
  title={Cluster characters for 2-Calabi--Yau triangulated categories},
  author={Palu, Yann},
  journal={Annales de l'Institut Fourier},
  volume={58},
  number={6},
  pages={2221--2248},
  year={2008}
}

@article{caldero2006triangulatedII,
  title={From triangulated categories to cluster algebras {II}},
  author={Caldero, Philippe and Keller, Bernhard},
  journal={Annales Scientifiques de l'\'Ecole Normale Sup\'erieure},
  volume={39},
  number={6},
  pages={983--1009},
  year={2006}
}

@article{zhou2014cluster,
  title={Cluster algebras arising from cluster tubes},
  author={Zhou, Yu and Zhu, Bin},
  journal={Journal of the London Mathematical Society},
  volume={89},
  number={3},
  pages={703--723},
  year={2014},
  publisher={Oxford University Press}
}

@article{reading2026theta,
  title={Theta functions in acyclic affine type},
  author={Reading, Nathan and Stella, Salvatore},
  journal={arXiv preprint arXiv:2603.23429},
  year={2026}
}

@article{fu2020cluster,
  title={Cluster algebras arising from cluster tubes II: The Caldero-Chapoton map},
  author={Fu, Changjian and Geng, Shengfei and Liu, Pin},
  journal={Journal of Algebra},
  volume={544},
  pages={228--261},
  year={2020},
  publisher={Elsevier}
}

@article{labardini2024gentle,
  title={Gentle Algebras Arising from Surfaces with Orbifold Points of Order 3, Part I: Scattering Diagrams},
  author={Labardini-Fragoso, Daniel and Mou, Lang},
  journal={Algebras and Representation Theory},
  volume={27},
  number={1},
  pages={679--722},
  year={2024},
  publisher={Springer}
}

@book{Ringel1984Tame,
  author    = {Ringel, Claus Michael},
  title     = {Tame Algebras and Integral Quadratic Forms},
  series    = {Lecture Notes in Mathematics},
  volume    = {1099},
  publisher = {Springer},
  year      = {1984}
}

@article{keller2008cluster,
  title={Cluster algebras, quiver representations and triangulated categories},
  author={Keller, Bernhard},
  journal={arXiv preprint arXiv:0807.1960},
  year={2008}
}

@article{gross2018canonical,
  title={Canonical bases for cluster algebras},
  author={Gross, Mark and Hacking, Paul and Keel, Sean and Kontsevich, Maxim},
  journal={Journal of the American Mathematical Society},
  volume={31},
  number={2},
  pages={497--608},
  year={2018}
}

@article{buan2010cluster,
  title={Cluster structures from 2-Calabi--Yau categories with loops},
  author={Buan, Aslak Bakke and Marsh, Bethany R and Vatne, Dagfinn F},
  journal={Mathematische Zeitschrift},
  volume={265},
  number={4},
  pages={951--970},
  year={2010},
  publisher={Springer}
}

@article{fomin2002cluster,
  author  = {Fomin, Sergey and Zelevinsky, Andrei},
  title   = {Cluster Algebras {I}: Foundations},
  journal = {Journal of the American Mathematical Society},
  volume  = {15},
  number  = {2},
  pages   = {497--529},
  year    = {2002}
}

@article{ding2013integral,
  title={Integral bases of cluster algebras and representations of tame quivers},
  author={Ding, Ming and Xiao, Jie and Xu, Fan},
  journal={Algebras and Representation Theory},
  volume={16},
  number={2},
  pages={491--525},
  year={2013}
}
\end{document}